\documentclass{jonasart}

\usepackage{pgfplots} \pgfplotsset{compat=1.18}

\newcommand{\A}{\mathcal{A}}
\newcommand{\af}{\alpha}
\newcommand{\bl}{\bullet}
\newcommand{\bt}{\beta}
\newcommand{\C}{\mathcal{C}}
\newcommand{\D}{\mathcal{D}}
\newcommand{\dl}{\delta}
\newcommand{\G}{\Gamma}
\newcommand{\impl}{\Rightarrow}
\newcommand{\J}{\mathcal{J}}

\newcommand{\m}{{}^{-1}}
\newcommand{\sg}{\sigma}
\newcommand{\sst}{\subseteq}
\newcommand{\vf}{\varphi}

\DeclareMathOperator{\Ann}{Ann}
\DeclareMathOperator{\Aut}{Aut}
\DeclareMathOperator{\id}{id}
\DeclareMathOperator{\spn}{span}

\DeclareFontFamily{U}{mathb}{}
\DeclareFontShape{U}{mathb}{m}{n}{<-5.5> mathb5 <5.5-6.5> mathb6 
	<6.5-7.5> mathb7 <7.5-8.5> mathb8 <8.5-9.5> mathb9 <9.5-11> mathb10 
	<11-> mathb12}{}
\DeclareSymbolFont{mathb}{U}{mathb}{m}{n}
\DeclareFontSubstitution{U}{mathb}{m}{n}
\DeclareMathSymbol{\blackdiamond}{\mathbin}{mathb}{"0C}

\title{Totally compatible structures on the radical of an incidence algebra}

\authors{Mykola Khrypchenko}

\authorinfo{Federal University of Santa Catarina, Brazil}{nskhripchenko@gmail.com}

\abstract{We describe totally compatible structures on the Jacobson radical of the incidence algebra of a finite poset over a field. We show that such structures are in general non-proper.}



\acknowledgments{The author thanks the referees for the careful reading of the manuscript and for the valuable comments.}

\editinfo{January 1, 2026}{March 26, 2026}{Samuel Lopes}

\VOLUME{1}
\ISSUE{1}
\NUMBER{2}
\DOI{https://doi.org/10.46298/jonas.17219}

\begin{document}

    \clearpage
    \section*{Introduction}

	Two algebraic structures from a~variety $\A$ defined on a~same vector space are said to be \textit{compatible} if their sum (equivalently, any linear combination) belongs to $\A$. Originated in the 70's~\cite{Magri78}, compatibility of Lie~\cite{Golubchik-Sokolov02,Golubchik-Sokolov05} associative~\cite{Carinena-Grabowski-Marmo2000,Odesskii-Sokolov06} and Poisson~\cite{Bolsinov91,Bolsinov-Borisov02,Reyman-Semenov89} algebras has been actively studied in the context of mathematical physics. Recently, it has attracted the interest of scientists from other areas of mathematics who expanded the investigation of compatibility to the varieties of $L_\infty$-algebras~\cite{Das22}, Leibniz algebras~\cite{Makhlouf-Saha23}, pre-Lie algebras~\cite{Abdelwahab-Kaygorodov-Makhlouf24}, anti-pre-Lie algebras~\cite{Normatov24}, left-symmetric algebras~\cite{Wu25}, Hom-Lie algebras~\cite{Das23CompHomLie}, Hom-Lie triple systems~\cite{Teng-Long-Zhang-Lin23} and so on. In particular,~\cite{Abdelwahab-Kaygorodov-Makhlouf24,Ladra-Cunha-Lopes24,Makhlouf-Saha23} provide algebraic classifications of low-dimensional compatible algebras.

	Given a~structure $\mu$ from $\A$, it is natural to ask: can one explicitly describe all the structures from $\A$ (on the same vector space) that are compatible with $\mu$? This question has been previously studied~\cite{Carinena-Grabowski-Marmo2000,Odesskii-Sokolov06} for the variety of associative algebras. Recall that the compatibility of two \textit{associative} bilinear products $\cdot_1$ and $\cdot_2$ on a~vector space $V$ is equivalent to the following identity on $V$:
	\begin{align}\label{compatible-ass-identity}
		(a\cdot_1 b)\cdot_2 c + (a\cdot_2 b)\cdot_1 c = a\cdot_1 (b\cdot_2 c) + a\cdot_2 (b\cdot_1 c).
	\end{align}
	In general, fixing $\cdot_1$, it is technically difficult to give a~complete description of $\cdot_2$ satisfying~\eqref{compatible-ass-identity}, although there are some natural classes of such products associated to $\cdot_1$.

	In this context, it seems to be more reasonable to study particular cases of \eqref{compatible-ass-identity} that are based on equalities of pairs of monomials. More precisely, let $\sg$ be a~permutation of $\{1,2\}$. If the $i$-th monomial of the left-hand side of \eqref{compatible-ass-identity} equals the $\sg(i)$-th monomial of the right-hand side of \eqref{compatible-ass-identity} for all $i\in\{1,2\}$, then $\cdot_1$ and $\cdot_2$ are said to be \textit{$\sg$-matching}~\cite{zgg23} (in the case $\sg=\id$ the triple $(V,\cdot_1,\cdot_2 )$ is called a~\textit{matching dialgebra}~\cite{zbg12}, also known as \textit{$As^{(2)}$-algebra}~\cite{Zinbiel}). If all the $4$ monomials in \eqref{compatible-ass-identity} are equal, then $\cdot_1$ and $\cdot_2$ are \textit{totally compatible} (see~\cite{Zhang13,zbg13}). In~\cite{Khr2024}, inspired by~\cite{Strohmayer08}, we introduced another similar notion, called \textit{interchangeability}, which means that the operations $\cdot_1$ and $\cdot_2$ can be permuted in each of the monomials in \eqref{compatible-ass-identity} without changing the brackets. If $A$ is \textit{unital} (or in certain sense close to being unital), then it is easy to see that $\sg$-matching, interchangeable and totally compatible structures on $A$ are simple modifications of the original associative product on $A$ (see~\cite{Khr2024}).

	The situation becomes more interesting in the context of \textit{nilpotent} associative algebras, where new classes of compatible structures arise. In~\cite{Khr2024b} we characterized $\sg$-matching, interchangeable and totally compatible structures on the strictly upper triangular matrix algebra $UT_n (K)$, $n\ge 3$, which is a~classical example of a~nilpotent associative algebra. Notice that $UT_n (K)$ is the \textit{Jacobson radical} of the algebra $T_n (K)$ of \textit{all} upper triangular $n\times n$ matrices over $K$, which is in turn a~particular case of the \textit{incidence algebra} $I(X,K)$ of a~\textit{finite} poset $X$ over $K$. So, our next goal was to generalize the results of~\cite{Khr2024b} to the Jacobson radical $J(I(X,K))$ of $I(X,K)$. However, having classified the $\id$-matching structures on $J(I(X,K))$ in the first draft of this manuscript, we gave up on the idea to proceed to $(12)$-matching and interchangeable structures, since the classification turned out to be too complex. We thus focused only on \textit{totally compatible} structures that are the main object of study in this article.

	In Section~\ref{sec-prelim} we collect all the needed definitions and preliminary results about totally compatible structures and incidence algebras to make the paper self-contained.

	In Section~\ref{sec-gen-prop} we prove a~couple of short general results whose particular cases will be used below.

	In Section~\ref{sec-proper-general} we introduce \textit{annihilator-valued} structures and study their relationship with structures $*_\vf$ determined by centroid elements (known from~\cite{Khr2024}). As a~result, we arrive at the definition of a~\textit{proper} totally compatible structure $*$, which essentially means that, up to an annihilator-valued component, $*$ is determined by a~centroid element. The motivation comes from~\cite{Khr2024,Khr2024b}, where most of the algebras had only totally compatible structures of this type. In Proposition~\ref{proper-tot-comp-closed-by-iso-and-antiiso}, we prove that the class of proper totally compatible structures is closed under isomorphisms and antiisomorphisms.

	Having prepared all the necessary general background, in Section~\ref{sec-prop-tot-comp} we focus on the algebra $J(I(X,K))$ and first describe in Proposition~\ref{descr-centroid-J} its centroid, which is itself an interesting object of study. It is then used in Corollary~\ref{proper-tot-comp-explicit} to give a~description of proper totally compatible structures on $J(I(X,K))$.

	Finally, in Section~\ref{sec-tot-comp} we prove our main result Theorem~\ref{tot-comp-exactly-bl+lin-comb-*_C}, which gives a~complete description of totally compatible structures on $J(I(X,K))$. Annihilator-valued structures again play an important role in this case, but, instead of $*_\vf$, we obtain a~more general class of structures that behave locally as $*_\vf$. In Corollary~\ref{suff-cond-all-tot-comp-are-proper}, we find a~sufficient condition on $X$ under which all the totally compatible structures on $J(I(X,K))$ are proper. Example~\ref{Y-poset-non-proper} shows that non-proper totally compatible structures on $J(I(X,K))$ may indeed exist, and Example~\ref{Y-poset-non-proper} demonstrates that the condition from Corollary~\ref{suff-cond-all-tot-comp-are-proper} is not necessary. We then give Example~\ref{Y-poset-ann-val} that illustrates the situation when all the totally compatible structures on $J(I(X,K))$ are annihilator-valued. In fact, this happens if and only if $X$ has length at most $2$, as proved in Proposition~\ref{all-tot-comp-ann-val<=>l(X)<=2}. We leave as an open problem to describe finite posets $X$ such that all the totally compatible structures on $J(I(X,K))$ are proper.

	\section{Definitions and preliminaries}\label{sec-prelim}

	Throughout the paper $K$ will be an arbitrary field. All the algebras and vector spaces will be over $K$ and all the products will be $K$-bilinear. As usual, the symbol $\circ$ will be reserved for the composition of two maps.

	Let $(A,\cdot)$ be a~(not necessarily associative) algebra. The \textit{annihilator} of $A$ is the ideal
	\begin{align*}
		\Ann(A,\cdot)=\{a\in A: a\cdot b=b\cdot a=0\text{ for all }b\in A\}.
	\end{align*}
	We will often write simply $\Ann(A)$, when it is clear what product is meant. Given subsets $B,C\sst A$, we denote
	\begin{align*}
		B\cdot C=\spn_K\{b\cdot c: b\in B,c\in C\}.
	\end{align*}
	We may write $B^2$ for $B\cdot B$, when the product on $A$ is clear from the context. The \textit{centroid} of $A$, denoted $\G(A)$, is the space of linear maps $\vf:A\to A$ such that
	\begin{align}\label{vf(xy)=vf(x).y=x.vf(y)}
		\vf(a\cdot b)=a\cdot\vf(b)=\vf(a)\cdot b
	\end{align}
	for all $a,b\in A$.

	\subsection{Matching compatibilities and interchangeability}

	Generalizing~\cite[Definition 1.1]{Khr2024} (inspired by~\cite{zgg23}) we say that two (not necessarily associative!) products $\cdot_1$ and $\cdot_2$ on a~same vector space $V$ are
	\begin{enumerate}
		\item \textit{$\sg$-matching} (where $\sg\in S_2 =\{\id,(12)\}$), if
		\begin{align}
			(a\cdot_1 b)\cdot_2 c = a\cdot_{\sg(1)}(b\cdot_{\sg(2)}c) \text{ and } (a\cdot_2 b)\cdot_1
			c = a\cdot_{\sg(2)}(b\cdot_{\sg(1)}c);
		\end{align}
		\item \textit{interchangeable} if
		\begin{align*}
			(a\cdot_1 b)\cdot_2 c = (a\cdot_2 b)\cdot_1 c \text{ and } a\cdot_1 (b\cdot_2 c) = a\cdot_2 (b\cdot_1 c);
		\end{align*}
		\item \textit{totally compatible}, if
		\begin{align}
			(a\cdot_1 b)\cdot_2 c = (a\cdot_2 b)\cdot_1 c = a\cdot_1 (b\cdot_2 c) = a\cdot_2 (b\cdot_1 c),
		\end{align}
	\end{enumerate}
	for all $a,b,c\in V$. Observe that $\sg$-matching \textit{associative} products are compatible in the sense of \eqref{compatible-ass-identity}, while interchangeable ones need not be compatible in general. Furthermore, notice that $\cdot_1$ and $\cdot_2$ are totally compatible if and only if $\cdot_1$ and $\cdot_2$ are $\sg$-matching \textit{for all} $\sg\in S_2$ if and only if $\cdot_1$ and $\cdot_2$ are interchangeable and $\sg$-matching \textit{for some} $\sg\in S_2$.

	Given an \textit{associative} algebra $(A,\cdot)$, by a~\textit{$\sg$-matching} (resp. \textit{interchangeable} or \textit{totally compatible}) \textit{structure} on $(A,\cdot)$ we mean an \textit{associative} product $*$ on $A$ such that $\cdot$ and $*$ are $\sg$-matching (resp. interchangeable or totally compatible). It is obvious that $\cdot$ is a~totally compatible structure on $(A,\cdot)$. It is known that any \textit{mutation}~\cite{ElduqueMyung} $\cdot_x$ of $\cdot$ by $x\in A$, i.e. $a\cdot_x b=a\cdot x\cdot b$, is an $\id$-matching structure on $(A,\cdot)$. If, moreover, $x$ is central, then $\cdot_x$ is a~totally compatible structure on $(A,\cdot)$.

	As we saw in~\cite[Propositions 2.2 and 2.6]{Khr2024}, for a~\textit{unital} algebra $A$, any $\id$-matching (resp. $(12)$-matching or interchangeable) structure on $(A,\cdot)$ is a~mutation of $\cdot$ by $x\in A$ (resp. by a~central $x\in A$). Consequently, any $(12)$-matching or interchangeable structure on a~unital algebra $(A,\cdot)$ is totally compatible. If $A$ is non-unital (in particular, if $A$ is \textit{nilpotent}), all these compatibilities may be different (see~\cite[Examples 1.4 and 1.5]{Khr2024}).

	\subsection{Isomorphic and antiisomorphic structures}
	Let $*_1$ and $*_2$ be two $\sg$-matching, interchangeable or totally compatible structures on an associative algebra $(A,\cdot)$. Following~\cite{Khr2024,Khr2024b}, we say that an automorphism (resp. antiautomorphism) $\phi$ of $(A,\cdot)$ is an \textit{isomorphism} (resp. \textit{antiisomorphism}) between $*_1$ and $*_2$, if
	\begin{align*}
		\phi(a*_1 b)=\phi(a)*_2\phi(b)\ (\text{resp. }\phi(a*_1 b)=\phi(b)*_2\phi(a))
	\end{align*}
	for all $a,b\in A$. The structures $*_1$ and $*_2$ are \textit{isomorphic} (resp. \textit{antiisomorphic}), if there is an isomorphism (resp. antiisomorphism) between them. This means that $(A,\cdot,*_1 )$ and $(A,\cdot,*_2 )$ are isomorphic (resp. antiisomorphic) as algebras with two multiplications. Given a~$\sg$-matching, interchangeable or totally compatible structure $*$ on $(A,\cdot)$ and an automorphism (resp. antiautomorphism) $\phi$ of $(A,\cdot)$, the product $\star$ on $A$ defined by
	\begin{align}
		a\star b=\phi(\phi\m(a)*\phi\m(b))\ (\text{resp. }a\star b=\phi(\phi\m(b)*\phi\m(a)))
	\end{align}
	for all $a,b\in A$, is a~structure on $(A,\cdot)$ of the same kind as $*$, and it is isomorphic (resp. antiisomorphic) to $*$ (see~\cite[Lemma 1.2]{Khr2024b}).

	\subsection{Posets}

	Let $(X,\le)$ be a~finite poset. As usual, we write $x<y$ to mean $x\le y$ and $x\ne y$. The binary relations $\ge$ and $>$ are inverse to $\le$ and $<$, respectively. Two elements $x,y\in X$ are said to be \textit{comparable}, if $x\le y$ or $x\ge y$. A nonempty subset $C\sst X$ is a~\textit{chain}, if any two elements $x,y\in C$ are comparable. A chain $C\sst X$ is \textit{maximal}, if there is no chain properly containing $C$. The \textit{length} of a~chain $C\sst X$ is $l(C):=|C|-1$. The \textit{length} of $X$ is
	\begin{align*}
		l(X):=\max\{l(C): C\text{ is a~chain in }X\}.
	\end{align*}
	Given $x\le y$, define
	\begin{align*}
		l(x,y):=l(\{z\in X: x\leq z\leq y\}).
	\end{align*}
	Denote by $\min(X)$ and $\max(X)$ the subsets of minimal and maximal elements of $X$, respectively. Finally, write
	\begin{align*}
		X^n_< =\{(x_1,\dots,x_n )\in X^n: x_1 <\dots<x_n\}.
	\end{align*}

	\subsection{Incidence algebras}

	Let $(X,\le)$ be a~finite poset and $K$ a~field. The \textit{incidence algebra} $I(X,K)$ of $X$ over $K$ (see~\cite{Rota64,SpDo}) is the associative $K$-algebra with basis $\{e_{xy}: x,y\in X, x\le y\}$ (called \textit{natural basis}) and bilinear multiplication
	\[
	e_{xy}\cdot e_{uv}=\dl_{yu}e_{xv}
	\]
	for all $x\leq y$ and $u\leq v$ in $X$. Here and below $\dl$ means the Kronecker delta. Given $f\in I(X,K)$, denote by $f(x,y)$ the coefficient of $e_{xy}$ in the linear combination
	\begin{align*}
		f=\sum_{x\le y}f(x,y)e_{xy}.
	\end{align*}
	The following formulas will be useful:
	\begin{align*}
		e_{xy}\cdot f=\sum_{v\ge y}f(y,v)e_{xv}\text{ and }f\cdot e_{xy}=\sum_{u\le x}f(u,x)e_{uy}.
	\end{align*}
	Recall from~\cite[Theorem 4.2.5]{SpDo} that the \textit{Jacobson radical} of $I(X,K)$ is
	\begin{align*}
		J(I(X,K))=\{f\in I(X,K):f(x,x)=0,\forall x\in X\}=\spn_K\{e_{xy}: x<y\}.
	\end{align*}
	The following facts are well-known and easy to prove:
	\begin{align*}
		\Ann(J(I(X,K)))&=\spn_K\{e_{xy}: \min(X)\ni x<y\in\max(X)\},\\
		J(I(X,K))\cdot J(I(X,K))&=\spn_K\{e_{xy}: l(x,y)>1\}.
	\end{align*}
	To shorten some formulas, we will write $\J$ for $J(I(X,K))$ below.

	\section{Some general properties of totally compatible structures}\label{sec-gen-prop}

	The following is clear.
	\begin{lemma}\label{lin-comb-tot-comp-are-tot-comp}
		Let $\{\cdot_i\}_{i\in I}$ be a~family of bilinear products on a~same vector space $V$. If $\{\cdot_i\}_{i\in I}$ are pairwise $\sg$-matching (resp., interchangeable or totally compatible), then any two finite linear combinations of $\{\cdot_i\}_{i\in I}$ are $\sg$-matching (resp., interchangeable or totally compatible).
	\end{lemma}
	The next result permits us to avoid unnecessary proofs of associativity in some cases.
	\begin{lemma}\label{*_1 +*_2 -tot-comp-str<=>*_2 -ass}
		Let $(A,\cdot)$ be an associative algebra and $*_1$, $*_2$ two bilinear products on $A$, where $*_1$ is associative. Assume that $\cdot$, $*_1$ and $*_2$ are pairwise $\sg$-matching (resp. totally compatible). Then $*_1 +*_2$ is a~$\sg$-matching (resp. totally compatible) structure on $(A,\cdot)$ if and only if $*_2$ is associative.
	\end{lemma}

	\begin{proof}
		It is enough to prove the result for $\sg$-matching structures.

		\textit{The ``only if'' part.} Assume that $*_1 +*_2$ is a~$\sg$-matching structure on $(A,\cdot)$. In particular, $*_1 +*_2$ is associative. By Lemma~\ref{lin-comb-tot-comp-are-tot-comp} the products $*_1 +*_2$ and $*_1$ are $\sg$-matching, whence they are compatible. It follows that $*_2 =(*_1 +*_2 )-*_1$ is associative.

		\textit{The ``if'' part.} Assume that $*_2$ is associative. Since $*_1$ and $*_2$ are $\sg$-matching and associative, their sum $*_1 +*_2$ is also associative. By Lemma~\ref{lin-comb-tot-comp-are-tot-comp} the products $*_1 +*_2$ and $\cdot$ are $\sg$-matching. Thus, $*_1 +*_2$ is a~$\sg$-matching structure on $(A,\cdot)$.
	\end{proof}

	The following property is a~particular case of total compatibility that will occur below.
	\begin{definition}
		Let $\cdot_1$ and $\cdot_2$ be two (not necessarily associative) products on a~same vector space $V$. We say that $\cdot_1$ and $\cdot_2$ are \textit{mutually annihilating} if
		\begin{align}
			(V\cdot_1 V)\cdot_2 V=(V\cdot_2 V)\cdot_1 V=V\cdot_1 (V\cdot_2 V) = V\cdot_2 (V\cdot_1 V)=\{0\}.
		\end{align}
	\end{definition}
	\section{Proper totally compatible structures}\label{sec-proper-general}

	Recall from~\cite[Definition 2.4]{Khr2024} the following.
	\begin{definition}
		Let $(A,\cdot)$ be a~(not necessarily associative) algebra and $\vf\in\G(A)$. The product $*_\vf$ on $A$ given by
		\begin{align}
			a*_\vf b:=\vf(a\cdot b)
		\end{align}
		is said to be \textit{determined} by $\vf$.
	\end{definition}

	\begin{remark}
		In particular, $\cdot$ is determined by $\id$.
	\end{remark}
	The next fact was proved as a~part of~\cite[Lemma 2.3]{Khr2024}.
	\begin{lemma}\label{*_v f-tot-comp-with-cdot}
		Let $(A,\cdot)$ be an associative algebra and $\vf\in\G(A)$. Then $*_\vf$
		is a~totally compatible structure on $(A,\cdot)$.
	\end{lemma}
	Let us introduce one more class of structures, which plays a~particularly important role in the case of nilpotent algebras.
	\begin{definition}
		Let $(A,\cdot)$ be a~(not necessarily associative) algebra. A bilinear product $\bl$ on $A$ is said to be \textit{annihilator-valued} (with respect to $\cdot$), if $\cdot$ and $\bl$ are mutually annihilating. An \textit{annihilator-valued structure} on an \textit{associative} algebra $(A,\cdot)$ is an \textit{associative} annihilator-valued bilinear product on $A$.
	\end{definition}

	\begin{remark}\label{ann-val-equiv-descr}
		A bilinear product $\bl$ on $A$ is annihilator-valued if and only if
		\begin{enumerate}
			\item\label{A-bl-A-sst-Ann(A)} $A\bl A\sst\Ann(A,\cdot)$;
			\item\label{(A-cdot-A)-bl-A=A-bl-(A-cdot-A)=0} $A\cdot A\sst\Ann(A,\bl)$.
		\end{enumerate}
	\end{remark}

	\begin{remark}\label{cdot-and-bl-mut-ann}
		Any annihilator-valued product $\bl$ on $A$ is totally compatible with $\cdot$.
	\end{remark}

	\begin{remark}\label{null-algebra-tot-comp}
		An annihilator-valued product $\bl$ on $A$ may be non-associative even if $(A,\cdot)$ is associative. For example, consider $A$ to be an algebra with the trivial product, i.e. $A\cdot A=\{0\}$, so that $\Ann(A,\cdot)=A$. Then any bilinear product $\bl$ on $A$ is annihilator-valued.
	\end{remark}

	\begin{remark}\label{bl-associative-if-A^2 -sst-Ann(A)}
		If $\Ann(A,\cdot)\sst A\cdot A$, then any annihilator-valued product $\bl$ on $A$ satisfies $A\bl A\sst\Ann(A,\bl)$, so $\bl$ is associative in this case, regardless of whether $(A,\cdot)$ is associative.
	\end{remark}
	Annihilator-valued products on $(A,\cdot)$ admit the following constructive description.
	\begin{lemma}
		Let $(A,\cdot)$ be a~(not necessarily associative) algebra and $V$ a~vector space complement of $A\cdot A$ in $A$. Then a~bilinear product $\bl$ on $A$ is annihilator-valued if and only if there exists a~bilinear map $\mu:V\times V\to\Ann(A,\cdot)$ such that
		\begin{align}\label{annihilator-valued-structure}
			(a_1 +a_2 )\bl(b_1 +b_2 )=\mu(a_2,b_2 ),
		\end{align}
		where $a_1,b_1\in A\cdot A$ and $a_2,b_2\in V$.
	\end{lemma}

	\begin{proof}
		\textit{The ``if'' part.} Let $\mu:V\times V\to\Ann(A,\cdot)$ be a~bilinear map and define $\bl$ by \eqref{annihilator-valued-structure}. It is clear that $\bl$ is bilinear and $A\bl A\sst\Ann(A,\cdot)$. Now, if $a\in A\cdot A$, then $a=a_1 +a_2$, where $a_1 =a\in A\cdot A$ and $a_2 =0\in V$. Hence, $a\bl b=\mu(0,b_2 )=0$ for any $b=b_1 +b_2\in A$ with $b_1\in A\cdot A$ and $b_2\in V$. Thus, $(A\cdot A)\bl A=\{0\}$. The proof that $A\bl (A\cdot A)=\{0\}$ is similar.

		\textit{The ``only if'' part.} Let $\bl$ be an annihilator-valued product on $A$. Given $a=a_1 +a_2$ and $b=b_1 +b_2$ with $a_1,b_1\in A\cdot A$ and $a_2,b_2\in V$, by Remark~\ref{ann-val-equiv-descr}\eqref{(A-cdot-A)-bl-A=A-bl-(A-cdot-A)=0} we have
		\begin{align*}
			(a_1 +a_2 )\bl(b_1 +b_2 )=a_2\bl b_2.
		\end{align*}
		So, we may define $\mu:V\times V\to\Ann(A,\cdot)$ by $\mu(u,v)=u\bl v$ for all $u,v\in V$. Then $\mu$ is bilinear and \eqref{annihilator-valued-structure} holds.
	\end{proof}

	\begin{lemma}\label{mixed-monomials-*_v f-bl-are-zero}
		Let $(A,\cdot)$ be a~(not necessarily associative) algebra, $\vf\in\G(A)$ and $\bl$ an annihilator-valued product on $A$. Then $*_\vf$ and $\bl$ are mutually annihilating, in particular, totally compatible.
	\end{lemma}

	\begin{proof}
		For all $a,b,c\in A$ we have
		\begin{align*}
			(a*_\vf b)\bl c=\vf(a\cdot b)\bl c=(a\cdot\vf(b))\bl c=0
		\end{align*}
		by Remark~\ref{ann-val-equiv-descr}\eqref{(A-cdot-A)-bl-A=A-bl-(A-cdot-A)=0} and
		\begin{align*}
			(a\bl b)*_\vf c=\vf((a\bl b)\cdot c)=\vf(0)=0
		\end{align*}
		by Remark~\ref{ann-val-equiv-descr}\eqref{A-bl-A-sst-Ann(A)}. Symmetrically, $a*_\vf (b\bl c) = a\bl(b*_\vf c)=0$.
	\end{proof}

	Hence, as a~consequence of Lemmas~\ref{*_v f-tot-comp-with-cdot},~\ref{mixed-monomials-*_v f-bl-are-zero},~\ref{*_1 +*_2 -tot-comp-str<=>*_2 -ass}, and Remark~\ref{cdot-and-bl-mut-ann}, we get.
	\begin{corollary}
		Let $(A,\cdot)$ be an associative algebra, $\vf\in\G(A)$ and $\bl$ an annihilator-valued product on $A$. Then $*_\vf+\bl$ is a~totally compatible structure on $(A,\cdot)$ if and only if $\bl$ is associative, i.e. $\bl$ is an annihilator-valued structure on $(A,\cdot)$.
	\end{corollary}

	\begin{definition}\label{defn-proper-str}
		Let $(A,\cdot)$ be an associative algebra. A totally compatible structure on $(A,\cdot)$ is said to be \textit{proper}, if it is of the form $*_\vf+\bl$, where $\vf\in\G(A)$ and $\bl$ is an annihilator-valued structure on $(A,\cdot)$.
	\end{definition}

	In most cases studied before in~\cite{Khr2024,Khr2024b} all the totally compatible structures turned out to be proper.
	\begin{example}
		Let $(A,\cdot)$ be an associative algebra and $*$ a~totally compatible structure on $(A,\cdot)$.
		\begin{enumerate}
			\item If $A$ is unital, then $*$ is determined by $\vf\in\G(A)\cong C(A)$, where $C(A)$ is the center of $A$. Hence, $*$ is proper (see~\cite[Proposition~2.6]{Khr2024}).
			\item If $A$ is an algebra with the trivial product, then $*$ is annihilator-valued, whence $*$ is proper (see Remark~\ref{null-algebra-tot-comp}).
			\item If $A$ is the semigroup algebra of a~rectangular band, then $*$ is determined by $\vf\in\G(A)$, whence $*$ is proper (see~\cite[Proposition~3.11]{Khr2024}).
			\item If $A$ has enough idempotents, then $*$ is determined by $\vf\in\G(A)$, whence $*$ is proper (see~\cite[Proposition 3.23]{Khr2024}).
			\item If $A$ is the free non-unital algebra over a~set $X$ with $|X|>1$, then $*$ is determined by $\vf\in\G(A)$, whence $*$ is proper (see~\cite[Proposition~4.6]{Khr2024}).
			\item If $A$ is the free non-unital commutative algebra over a~set $X$ with $|X|>1$, then $*$ is determined by $\vf\in\G(A)$, whence $*$ is proper (see~\cite[Proposition~4.12]{Khr2024}).
			\item If $A$ is the strictly upper triangular matrix algebra, then $*$ is proper (see~\cite[Theorem~5.1]{Khr2024b}, where any linear combination of the structures $\mathbf{T}_{i,j}^1$ is annihilator-valued and any scalar multiple of the structure $\mathbf{T}^2$ is determined by $\vf\in\G(A)$).
			\item If $A$ is the free non-unital algebra over a~set $X$ with $|X|=1$ (i.e., the non-unital polynomial algebra in one variable), then $*$ is not always proper (see~\cite[Remark~4.10]{Khr2024}).
		\end{enumerate}
	\end{example}

	\begin{proposition}\label{proper-tot-comp-closed-by-iso-and-antiiso}
		The class of proper totally compatible structures on an associative algebra $(A,\cdot)$ is closed under isomorphisms and antiisomorphisms.
	\end{proposition}

	\begin{proof}
		Let $*=*_\vf+\bl$, where $\vf\in\G(A)$ and $\bl$ is an annihilator-valued structure on $(A,\cdot)$.

		Any structure $\star$ isomorphic to $*$ is of the form $a\star b=\phi(\phi\m(a)*\phi\m(b))$ for some $\phi\in\Aut(A,\cdot)$. Then for all $a,b\in A$ we have
		\begin{align*}
			a\star b&=\phi(\phi\m(a)*\phi\m(b))=\phi(\phi\m(a)*_\vf\phi\m(b)+\phi\m(a)\bl\phi\m(b))\\
			&=\phi(\vf(\phi\m(a)\cdot\phi\m(b)))+\phi(\phi\m(a)\bl\phi\m(b))\\
			&=(\phi\circ\vf\circ\phi\m)(a\cdot b)+\phi(\phi\m(a)\bl\phi\m(b)).
		\end{align*}
		Now observe that $\phi\circ\vf\circ\phi\m\in\G(A)$, because
		\begin{align*}
			(\phi\circ\vf\circ\phi\m)(a\cdot b)&=\phi(\vf(\phi\m(a)\cdot \phi\m(b)))=\phi(\phi\m(a)\cdot\vf(\phi\m(b)))\\
			&=a\cdot\phi(\vf(\phi\m(b)))=a\cdot(\phi\circ\vf\circ\phi\m)(b)
		\end{align*}
		and similarly $(\phi\circ\vf\circ\phi\m)(a\cdot b)=(\phi\circ\vf\circ\phi\m)(a)\cdot b$. Furthermore, the bilinear product
        \[
            a\blackdiamond b=\phi(\phi\m(a)\bl\phi\m(b))
        \]
        on $A$ is clearly associative, because $\bl$ is, and satisfies $A\blackdiamond A\sst\Ann(A,\cdot)$ thanks to Remark~\ref{ann-val-equiv-descr}\eqref{A-bl-A-sst-Ann(A)} and $\phi(\Ann(A,\cdot))=\Ann(A,\cdot)$. One also has
		$(A\cdot A)\blackdiamond A=A\blackdiamond (A\cdot A)=\{0\}$ due to Remark~\ref{ann-val-equiv-descr}\eqref{(A-cdot-A)-bl-A=A-bl-(A-cdot-A)=0}, $\phi\m(A)=A$ and $\phi\m(A\cdot A)=A\cdot A$. Thus, $\blackdiamond$ is annihilator-valued, so that $\star$ is proper as being equal to $*_{\phi\circ\vf\circ\phi\m}+\blackdiamond$.

		If $\star$ is antiisomorphic to $*$, then it is given by $a\star b=\phi(\phi\m(b)*\phi\m(a))$ for some antiisomorphism $\phi$ of $(A,\cdot)$. Similarly to the isomorphic structure, we have
		\begin{align*}
			a\star b
			&=\phi(\vf(\phi\m(b)\cdot\phi\m(a)))+\phi(\phi\m(b)\bl\phi\m(a))\\
			&=(\phi\circ\vf\circ\phi\m)(a\cdot b)+\phi(\phi\m(b)\bl\phi\m(a)),
		\end{align*}
		where
		\begin{align*}
			(\phi\circ\vf\circ\phi\m)(a\cdot b)&=\phi(\vf(\phi\m(b)\cdot \phi\m(a)))=\phi(\vf(\phi\m(b))\cdot\phi\m(a))\\
			&=a\cdot\phi(\vf(\phi\m(b)))=a\cdot(\phi\circ\vf\circ\phi\m)(b)
		\end{align*}
		and $(\phi\circ\vf\circ\phi\m)(a\cdot b)=(\phi\circ\vf\circ\phi\m)(a)\cdot b$, so that $\phi\circ\vf\circ\phi\m\in\G(A)$. It is easily seen as above that $a\blackdiamond b=\phi(\phi\m(b)\bl\phi\m(a))$ is an annihilator-valued structure on $(A,\cdot)$. Thus, $*=*_{\phi\circ\vf\circ\phi\m}+\blackdiamond$ is proper.
	\end{proof}

	\section{Proper totally compatible structures on $(\J,\cdot)$}\label{sec-prop-tot-comp}

	\subsection{The centroid of $\J$}

	In order to characterize proper totally compatible structures on $(\J,\cdot)$, we need a~description of the centroid of $\J$. Let us introduce a~class of centroid elements that will be one of the ingredients of the future description.
	\begin{definition}
		Let $(A,\cdot)$ be an algebra (which is not necessarily associative). A linear map $\vf:A\to A$ is said to be \textit{annihilator-valued}, if it satisfies
		\begin{enumerate}
			\item $\vf(A)\sst \Ann(A)$;
			\item $\vf(A\cdot A)=\{0\}$.
		\end{enumerate}
	\end{definition}

	\begin{remark}\label{ann-val-centroid}
		Any annihilator-valued linear map $\vf:A\to A$ belongs to $\G(A)$, since all the products in \eqref{vf(xy)=vf(x).y=x.vf(y)} are zero.
	\end{remark}

	Another class of centroid elements has its origin in~\cite{KK8}. Recall from~\cite{KK8} that a~map $\sg:X^2_<\to K$ is said to be \textit{constant on chains} if
	\begin{align*}
		\sg(x,y)=\sg(u,v),\text{ whenever }x<y\text{ and }u<v\text{ belong to a~same chain in }X.
	\end{align*}

	\begin{definition}
		Let $\sg:X^2_<\to K$ be constant on chains. Define the linear map $\vf_\sg:\J\to \J$ by
		\begin{align*}
			\vf_\sg(e_{xy})=\sg(x,y)e_{xy},
		\end{align*}
		for all $x<y$.
	\end{definition}

	\begin{lemma}\label{vf_s g-in-G(J)}
		Let $\sg:X^2_<\to K$ be constant on chains. Then $\vf_\sg\in\G(\J)$.
	\end{lemma}

	\begin{proof}
		Let $x<y$ and $u<v$. If $v\ne x$, then
		\begin{align*}
			e_{uv}\cdot\vf_\sg(e_{xy})=\sg(x,y)e_{uv}\cdot e_{xy}=0=\vf_\sg(0)=\vf_\sg(e_{uv}\cdot e_{xy}).
		\end{align*}
		Otherwise, if $v=x$, then $x<y$ and $u<y$ belong to the same chain, $u<x<y$. Hence, $\sg(x,y)=\sg(u,y)$, so
		\begin{align*}
			e_{uv}\cdot\vf_\sg(e_{xy})=\sg(x,y)e_{uv}\cdot e_{xy}=\sg(x,y)e_{uy}=\sg(u,y)e_{uy}=\vf_\sg(e_{uy})=\vf_\sg(e_{uv}\cdot e_{xy}).
		\end{align*}
		Similarly, one proves that $\vf_\sg(e_{xy})e_{uv}=\vf_\sg(e_{xy}\cdot e_{uv})$. Thus, $\vf_\sg\in\G(\J)$.
	\end{proof}

	\begin{lemma}\label{vf(e_x y)-in-spn-e_x y+Ann(J)}
		Let $\vf\in\G(\J)$ and $x<y$.
		\begin{enumerate}
			\item\label{vf(e_x y)-in-span(e_x y)} If $l(x,y)>1$, then $\vf(e_{xy})\in\spn_K\{e_{xy}\}$.
			\item\label{vf(e_x y)-in-span(e_x y)+Ann} If $l(x,y)=1$, then $\vf(e_{xy})\in\spn_K\{e_{xy}\}+\Ann(\J)$.
		\end{enumerate}
	\end{lemma}

	\begin{proof}
		We first prove part~\eqref{vf(e_x y)-in-span(e_x y)}. If $l(x,y)>1$, then there is $x<z<y$, so
		\begin{align*}
			\vf(e_{xy})=\vf(e_{xz}\cdot e_{zy})=e_{xz}\cdot\vf(e_{zy})\in\spn_K\{e_{xv}:z<v\}.
		\end{align*}
		Similarly,
		\begin{align*}
			\vf(e_{xy})=\vf(e_{xz})\cdot e_{zy}\in\spn_K\{e_{uy}:u<z\}.
		\end{align*}
		Since
		\begin{align*}
			\spn_K\{e_{xv}:z<v\}\cap\spn_K\{e_{uy}:u<z\}=\spn_K\{e_{xy}\},
		\end{align*}
		the proof of part~\eqref{vf(e_x y)-in-span(e_x y)} is complete.

		Now we prove part~\eqref{vf(e_x y)-in-span(e_x y)+Ann}. Let $u<v$ with $(u,v)\ne(x,y)$. Assume first that $u\not\in\min(X)$ and choose $a<u$. Then
		\begin{align}\label{vf(e_x y)(u_v )=vf(e_a u.e_x y)(a_v )}
			\vf(e_{xy})(u,v)=(e_{au}\cdot\vf(e_{xy}))(a,v)=\vf(e_{au}\cdot e_{xy})(a,v).
		\end{align}
		We have two cases.

		\textit{Case 1:} $u\ne x$. Then $e_{au}\cdot e_{xy}=0$, so $\vf(e_{xy})(u,v)=0$ by \eqref{vf(e_x y)(u_v )=vf(e_a u.e_x y)(a_v )}.

		\textit{Case 2:} $u=x$ and $v\ne y$. Then $e_{au}\cdot e_{xy}=e_{ay}$, where $l(a,y)>1$. By part~\eqref{vf(e_x y)-in-span(e_x y)} we have $\vf(e_{ay})\in\spn_K\{e_{ay}\}$. Since $v\ne y$, then $\vf(e_{ay})(a,v)=0$, so $\vf(e_{xy})(u,v)=0$ in view of \eqref{vf(e_x y)(u_v )=vf(e_a u.e_x y)(a_v )}. One similarly shows that $\vf(e_{xy})(u,v)=0$, whenever $v\not\in\max(X)$. This completes the proof of part~\eqref{vf(e_x y)-in-span(e_x y)+Ann}.
	\end{proof}

	\begin{lemma}\label{from-vf-to-sg}
		Let $\vf\in\G(\J)$. Then the associated map $\sg:X^2_<\to K$, given by
		\begin{align*}
			\sg(x,y)=\vf(e_{xy})(x,y),
		\end{align*}
		for all $x<y$, is constant on chains.
	\end{lemma}

	\begin{proof}
		Let $x<y$ and $C$ be a~chain containing $x$ and $y$. Denote the least element of $C$ by $a$ and the greatest element of $C$ by $b$. We first show that $\sg(x,y)=\sg(a,y)$. It suffices to consider the case $a\ne x$. Then $a<x$, so
		\begin{align*}
			\vf(e_{xy})(x,y)=(e_{ax}\cdot\vf(e_{xy}))(a,y)=\vf(e_{ax}\cdot e_{xy})(a,y)=\vf(e_{ay})(a,y),
		\end{align*}
		whence $\sg(x,y)=\sg(a,y)$. Similarly, one proves that $\sg(a,y)=\sg(a,b)$. This implies that $\sg(x,y)=\sg(a,b)$.
	\end{proof}

	\begin{proposition}\label{descr-centroid-J}
		The elements of $\G(\J)$ are exactly the maps of the form $\vf_\sg+\eta$, where $\sg:X^2_<\to K$ is constant on chains and $\eta$ is annihilator-valued.
	\end{proposition}

	\begin{proof}
		Let $\vf\in\G(\J)$. Then there is $\sg:X^2_<\to K$ from Lemma~\ref{from-vf-to-sg} which is constant on chains. Define $\eta=\vf-\vf_\sg$. In view of Lemma~\ref{vf(e_x y)-in-spn-e_x y+Ann(J)} we have $\eta(e_{xy})\in\Ann(\J)$ and $\eta(e_{xy})=0$, whenever $l(x,y)>1$. The latter means that $\eta(\J\cdot \J)=\{0\}$. Thus, $\eta$ is annihilator-valued.

		Conversely, $\vf_\sg\in\G(\J)$ by Lemma~\ref{vf_s g-in-G(J)} and $\eta\in\G(\J)$ by Remark~\ref{ann-val-centroid}. Therefore, $\vf_\sg+\eta\in\G(\J)$.
	\end{proof}

	\subsection{The description of proper totally compatible structures}
	\begin{corollary}\label{proper-tot-comp-explicit}
		Proper totally compatible structures on $(\J,\cdot)$ are exactly bilinear products $*$ of the form
		\begin{align}\label{e_x y*e_u v-proper}
			e_{xy}*e_{uv}=\sg(x,v)\dl_{yu}e_{xv}+e_{xy}\bl e_{uv},
		\end{align}
		where $\sg:X^2_<\to K$ is constant on chains and $\bl$ is an annihilator-valued structure on $(\J,\cdot)$.
	\end{corollary}

	\begin{proof}
		This follows from Definition~\ref{defn-proper-str} and Proposition~\ref{descr-centroid-J}, because we have
		\begin{align*}
			e_{xy}*_\vf e_{uv}=\vf_\sg(e_{xy}\cdot e_{uv})+\eta(e_{xy}\cdot e_{uv})=\vf_\sg(\dl_{yu}e_{xv})=\sg(x,v)\dl_{yu}e_{xv},
		\end{align*}
        for $\vf=\vf_\sg+\eta\in\G(\J)$, as in Proposition~\ref{descr-centroid-J}.
	\end{proof}

	\begin{remark}\label{proper-str-as-bl+lin-comb}
		Let $\sim$ be the equivalence relation on $X^2_<$ generated by the pairs
        \[
            ((x,y),(u,v))\in X^2_<\times X^2_<,
        \]
        such that there exists a~chain in $X$ containing both $x<y$ and $u<v$. Then a~map $X^2_<\to K$ is constant on chains if and only if it is constant on $\sim$-classes of $X^2_<$.

		Taking a~$\sim$-class $\C$ of $X^2_<$, one can define $\sg_\C:X^2_<\to K$ by
		\begin{align*}
			\sg_\C(x,y)=
			\begin{cases}
				1, & (x,y)\in\C,\\
				0, & (x,y)\not\in\C.
			\end{cases}
		\end{align*}
		Then $\sg_\C$ is constant on chains, and the corresponding totally compatible structure $*_\C$ on $(\J,\cdot)$ determined by $\vf_{\sg_\C}\in\G(\J)$ is given by
		\begin{align*}
			e_{xy}*_\C e_{uv}=
			\begin{cases}
				e_{xv}, & y=u\text{ and }(x,v)\in\C,\\
				0, & \text{otherwise}.
			\end{cases}
		\end{align*}
		It is therefore easily seen by Corollary~\ref{proper-tot-comp-explicit} that the totally compatible structures on $(\J,\cdot)$ are exactly bilinear products $*$ that are sums of an annihilator-valued structure $\bl$ on $(\J,\cdot)$ and a~linear combination of the structures $*_\C$, where $\C$ runs through the set $X^2_< /{\sim}$ of $\sim$-classes of $X^2_<$.
	\end{remark}

	\section{Totally compatible structures on $(\J,\cdot)$}\label{sec-tot-comp}

	\subsection{A class of totally compatible structures}

	We first define a~combinatorial notion that will play a~role in the main result.
	\begin{definition}
		Let $(x,y,z),(u,v,w)\in X^3_<$. The triples $(x,y,z)$ and $(u,v,w)$ are said to be \textit{chained} if there is a~chain in $X$ containing both $x<y<z$ and $u<v<w$. Let $\approx$ be the equivalence relation generated by all the pairs of chained triples in $X^3_<$ and denote by $X^3_< /{\approx}$ the set of $\approx$-equivalence classes of $X^3_<$.
	\end{definition}

	\begin{remark}
		If $(x,y,z)$ and $(u,v,w)$ are chained, then $(x,z)\sim (u,w)$. Consequently, for any $(x,y,z),(u,v,w)\in X^3_<$ we have
		\begin{align}\label{(x_y_z )-approx(u_v_w )=>(x_z )-sim-(u_w )}
			(x,y,z)\approx(u,v,w)\impl(x,z)\sim(u,w).
		\end{align}
		The converse of \eqref{(x_y_z )-approx(u_v_w )=>(x_z )-sim-(u_w )} is not true in general (see Examples~\ref{Y-poset-non-proper}, \ref{Y-poset-proper}, \ref{Y-poset-ann-val} below).
	\end{remark}

	\begin{lemma}\label{*_C -tot-comp-with-cdot}
		Let $\C\in X^3_< /{\approx}$. Then the bilinear product $*_\C$ given by
		\begin{align}\label{e_x y*_C e_u v=e_x z-or-0}
			e_{xy}*_\C e_{uv}=
			\begin{cases}
				e_{xv}, & y=u\text{ and }(x,y,v)\in\C,\\
				0, & \text{otherwise},
			\end{cases}
		\end{align}
		is a~totally compatible structure on $(\J,\cdot)$.
	\end{lemma}

	\begin{proof}
		Let $x<y$, $z<u$ and $v<w$. We are going to show that
		\begin{align}\label{(e_x y.e_z u)*_C e_v w-assoc}
			(e_{xy}\cdot e_{zu})*_\C e_{vw}=(e_{xy}*_\C e_{zu})\cdot e_{vw}=e_{xy}\cdot (e_{zu}*_\C e_{vw})=e_{xy}*_\C (e_{zu}\cdot e_{vw}).
		\end{align}
		\textit{Case 1:} $y\ne z$. Then $e_{xy}\cdot e_{zu}=e_{xy}*_\C e_{zu}=0$, so the first two products of \eqref{(e_x y.e_z u)*_C e_v w-assoc} are zero. By \eqref{e_x y*_C e_u v=e_x z-or-0} the product $e_{zu}*_\C e_{vw}$ is either $e_{zw}$ or $0$. In any case, $e_{xy}\cdot (e_{zu}*_\C e_{vw})=0$. Similarly, $e_{zu}\cdot e_{vw}$ is either $e_{zw}$ or $0$. In any case, $e_{xy}*_\C (e_{zu}\cdot e_{vw})=0$ by \eqref{e_x y*_C e_u v=e_x z-or-0}.

		\textit{Case 2:} $y=z$. Then $(e_{xy}\cdot e_{zu})*_\C e_{vw}=e_{xu}*_\C e_{vw}$. We have the following two subcases.

		\textit{Case 2.1:} $u\ne v$. Then $e_{xu}*_\C e_{vw}=e_{zu}\cdot e_{vw}=e_{zu}*_\C e_{vw}=0$, so all the products of \eqref{(e_x y.e_z u)*_C e_v w-assoc}, except possibly for the second one, are zero. Since $e_{xy}*_\C e_{zu}$ is either $e_{xu}$ or $0$, we have $(e_{xy}*_\C e_{zu})\cdot e_{vw}=0$ as well.

		\textit{Case 2.2:} $u=v$. Observe that
		\begin{align}\label{(xyu)-(xuw)-(yuw)-(xyw)}
			(x,y,u)\approx(x,u,w)\approx(y,u,w)\approx(x,y,w),
		\end{align}
		because $x<y<u<w$. We have the following two subcases.

		\textit{Case 2.2.1:} $(x,u,w)\not\in\C$. Then $(x,y,u),(y,u,w),(x,y,w)\not\in\C$. By \eqref{e_x y*_C e_u v=e_x z-or-0} we have $e_{xu}*_\C e_{vw}=0$, so that $(e_{xy}\cdot e_{zu})*_\C e_{vw}=0$. Since $e_{xy}*_\C e_{zu}=e_{zu}*_\C e_{vw}=0$, then $(e_{xy}*_\C e_{zu})\cdot e_{vw}=e_{xy}\cdot (e_{zu}*_\C e_{vw})=0$. Furthermore, $e_{xy}*_\C (e_{zu}\cdot e_{vw})=e_{xy}*_\C e_{zw}$, which is also zero by \eqref{e_x y*_C e_u v=e_x z-or-0}.

		\textit{Case 2.2.2:} $(x,u,w)\in\C$. Then $(x,y,u),(y,u,w),(x,y,w)\in\C$. In this case, we have $e_{xu}*_\C e_{vw}=e_{xw}$, so that $(e_{xy}\cdot e_{zu})*_\C e_{vw}=e_{xw}$ by \eqref{e_x y*_C e_u v=e_x z-or-0}. Now,         by \eqref{e_x y*_C e_u v=e_x z-or-0}:
        \begin{gather*}
            (e_{xy}*_\C e_{zu})\cdot e_{vw}
            =e_{xu}\cdot e_{vw}
            =e_{xw}, \\
            e_{xy}\cdot (e_{zu}*_\C e_{vw})
            =e_{xy}\cdot e_{zw}
            =e_{xw}, \\
            e_{xy}*_\C (e_{zu}\cdot e_{vw})
            =e_{xy}*_\C e_{zw}
            =e_{xw}.
        \end{gather*}
		Thus, the proof of \eqref{(e_x y.e_z u)*_C e_v w-assoc} is complete.

		For the associativity of $*$ let us show that
		\begin{align}\label{(e_x y*_C -e_z u)*_C e_v w-assoc}
			(e_{xy}*_\C e_{zu})*_\C e_{vw}=e_{xy}*_\C (e_{zu}*_\C e_{vw}).
		\end{align}
		The proof of \eqref{(e_x y*_C -e_z u)*_C e_v w-assoc} is similar to that of \eqref{(e_x y.e_z u)*_C e_v w-assoc}. As above, one sees that both of the monomials of \eqref{(e_x y*_C -e_z u)*_C e_v w-assoc} are zero, whenever $y\ne z$ or $u\ne v$. If $y=z$ and $u=v$, then $x<y<u<w$, so that \eqref{(xyu)-(xuw)-(yuw)-(xyw)} holds. Hence, both of the monomials of \eqref{(e_x y*_C -e_z u)*_C e_v w-assoc} are either zero (if $(x,y,u)\not\in\C$) or are equal to $e_{xw}$ (if $(x,y,u)\in\C$).
	\end{proof}

	\begin{lemma}\label{*_C -tot-comp-with-bl}
		Let $\C$ be a~$\approx$-class and $\bl$ an annihilator-valued product on $\J$. Then $*_\C$ and $\bl$ are mutually annihilating, in particular, totally compatible.
	\end{lemma}

	\begin{proof}
		It follows from $\J\bl \J\sst \Ann(\J,\cdot)$ and \eqref{e_x y*_C e_u v=e_x z-or-0} that $(\J\bl \J)*_\C \J=\J*_\C (\J\bl \J)=\{0\}$. Now, since $\J*_\C\J\sst\J\cdot\J$ and $\J\cdot\J\sst\Ann(\J,\bl)$, then $\J*_\C\J\sst\Ann(\J,\bl)$.
	\end{proof}

	\begin{lemma}\label{*_C -tot-comp-with-*_D}
		For any pair of distinct $\approx$-classes $\C$ and $\D$, the structures $*_\C$ and $*_\D$ are mutually annihilating, in particular, totally compatible.
	\end{lemma}

	\begin{proof}
		Let $x<y$, $z<u$ and $v<w$.

		\textit{Case 1:} $y\ne z$ or $u\ne v$. Then
		\begin{align}
			(e_{xy}*_\C e_{zu})*_\D e_{vw}=(e_{xy}*_\D e_{zu})*_\C e_{vw}&=e_{xy}*_\C (e_{zu}*_\D e_{vw})\notag\\
			&=e_{xy}*_\D (e_{zu}*_\C e_{vw})=0\label{(e_x y*_C e_z u)*_D e_v w=dots=0}
		\end{align}
		by \eqref{e_x y*_C e_u v=e_x z-or-0}.

		\textit{Case 2:} $y=z$ and $u=v$. Then $x<y<u<w$, so we have \eqref{(xyu)-(xuw)-(yuw)-(xyw)}.

		\textit{Case 2.1:} $(x,y,u)\in\C$. Then $(x,y,u),(x,u,w),(y,u,w),(x,y,w)\not\in\D$ by \eqref{(xyu)-(xuw)-(yuw)-(xyw)}, whence \eqref{(e_x y*_C e_z u)*_D e_v w=dots=0}.

		\textit{Case 2.2:} $(x,y,u)\in\D$. Then $(x,y,u),(x,u,w),(x,y,w),(y,u,w)\not\in\C$ by \eqref{(xyu)-(xuw)-(yuw)-(xyw)}, whence \eqref{(e_x y*_C e_z u)*_D e_v w=dots=0}.

		\textit{Case 2.3:} $(x,y,u)\not\in\C\sqcup\D$. Then \eqref{(e_x y*_C e_z u)*_D e_v w=dots=0} is immediate in view of \eqref{(xyu)-(xuw)-(yuw)-(xyw)}.
	\end{proof}

	\begin{proposition}\label{bl+lin-comb-*_C -tot-comp}
		The sum of an annihilator-valued structure $\bl$ on $(\J,\cdot)$ and a~linear combination of the structures $*_\C$, where $\C$ runs through the set of $\approx$-equivalence classes of $X^3_<$, is a~totally compatible structure on $(\J,\cdot)$.
	\end{proposition}

	\begin{proof}
		By Lemma~\ref{lin-comb-tot-comp-are-tot-comp}, Remark~\ref{cdot-and-bl-mut-ann} and Lemmas~\ref{*_C -tot-comp-with-*_D}, \ref{*_C -tot-comp-with-bl}, \ref{*_C -tot-comp-with-cdot}, any product of the form $\bl+\sum_{\C\in X^3_< /{\approx}}\af_\C *_\C$ is totally compatible with $\cdot$. Since $\bl$ and all the structures $*_\C$ are associative and pairwise totally compatible, then $\bl+\sum_{\C\in X^3_< /{\approx}}\af_\C *_\C$ is also associative.
	\end{proof}

	\subsection{The description of totally compatible structures}

	Recall the following fact.
	\begin{lemma}
		[Lemma 4.1 from~\cite{Khr2024b}]
		\label{a-cdot-b=0=>a*b-in-Ann(A)}
		Let $(A,\cdot)$ be a~(not necessarily associative) algebra and $*$ a~bilinear product on $A$ such that $*$ and $\cdot$ are interchangeable. For all $a,b\in A$ if $a\cdot b=0$, then $a*b\in\Ann(A,\cdot)$.
	\end{lemma}

	\begin{lemma}\label{e_x y*e_u v=0-for-y-ne-u-and-l(x_y )>1-or-l(u_v )>1}
		Let $*$ be a~bilinear product on $\J$ such that $*$ and $\cdot$ are $(12)$-matching. Given $x<y$ and $u<v$ with $y\ne u$, if $l(x,y)>1$ or $l(u,v)>1$, then $e_{xy}*e_{uv}=0$.
	\end{lemma}

	\begin{proof}
		Assume that $l(x,y)>1$ and choose $x<z<y$. Then
		\begin{align*}
			e_{xy}*e_{uv}=(e_{xz}\cdot e_{zy})*e_{uv}=e_{xz}*(e_{zy}\cdot e_{uv})=0.
		\end{align*}
		The case $l(u,v)>1$ is similar.
	\end{proof}

	\begin{lemma}\label{e_x y*e_y z-in-span-e_x z+Ann}
		Assume that $*$ is a~bilinear product on $\J$ such that $*$ and $\cdot$ are totally compatible. Let $x<y<z$.
		\begin{enumerate}
			\item\label{e_x y*e_y z-l(x_y )>1-or-l(y_z )>1} If $l(x,y)>1$ or $l(y,z)>1$, then $e_{xy}*e_{yz}\in\spn_K\{e_{xz}\}$.
			\item\label{e_x y*e_y z-l(x_y )=l(y_z )=1} If $l(x,y)=l(y,z)=1$, then $e_{xy}*e_{yz}\in\spn_K\{e_{xz}\}+\Ann(\J)$.
		\end{enumerate}
	\end{lemma}

	\begin{proof}
		We first prove part~\eqref{e_x y*e_y z-l(x_y )>1-or-l(y_z )>1}. Assume that $l(x,y)>1$ and choose $x<u<y$. Then
		\begin{align*}
			e_{xy}*e_{yz}=(e_{xu}\cdot e_{uy})*e_{yz}=e_{xu}\cdot(e_{uy}*e_{yz})\in\spn_K\{e_{xv}:u<v\}.
		\end{align*}
		On the other hand,
		\begin{align*}
			e_{xy}*e_{yz}=(e_{xu}\cdot e_{uy})*e_{yz}=(e_{xu}* e_{uy})\cdot e_{yz}\in\spn_K\{e_{wz}:w<y\}.
		\end{align*}
		Consequently,
		\begin{align*}
			e_{xy}*e_{yz}\in\spn_K\{e_{xv}:u<v\}\cap\spn_K\{e_{wz}:w<y\}=\spn_K\{e_{xz}\}.
		\end{align*}
		The case $l(y,z)>1$ is similar.

		Now we prove part~\eqref{e_x y*e_y z-l(x_y )=l(y_z )=1}. Let $l(x,y)=l(y,z)=1$. Take $u<v$ such that $(u,v)\ne (x,z)$. Assume first that $u\not\in\min(X)$ and choose $a<u$. Then
		\begin{align*}
			(e_{xy}*e_{yz})(u,v)=(e_{au}\cdot(e_{xy}*e_{yz}))(a,v).
		\end{align*}
		We have two cases.

		\textit{Case 1:} $u\ne x$. Then $e_{au}\cdot(e_{xy}*e_{yz})=(e_{au}\cdot e_{xy})*e_{yz}=0$.

		\textit{Case 2:} $u=x$ and $v\ne z$. By part~\eqref{e_x y*e_y z-l(x_y )>1-or-l(y_z )>1} we have
		\begin{align*}
			e_{au}\cdot(e_{xy}*e_{yz})=e_{ax}\cdot(e_{xy}*e_{yz})=e_{ax}*(e_{xy}\cdot e_{yz})=e_{ax}*e_{xz}\in\spn_K\{e_{az}\},
		\end{align*}
		because $l(x,z)>1$. Since $v\ne z$, we conclude that $(e_{au}\cdot(e_{xy}*e_{yz}))(a,v)=0$.

		Thus, $(e_{xy}*e_{yz})(u,v)=0$, whenever $(u,v)\ne (x,z)$ and $u\not\in\min(X)$. One similarly proves that $(e_{xy}*e_{yz})(u,v)=0$, whenever $(u,v)\ne (x,z)$ and $v\not\in\max(X)$.
	\end{proof}

	\begin{lemma}\label{(e_x y*e_y z)(x_z )=(e_u v*e_v w)(u_w )}
		Assume that $*$ is a~bilinear product on $\J$ such that $*$ and $\cdot$ are totally compatible. Given $x<y<z$ and $u<v<w$, if $(x,y,z)\approx(u,v,w)$, then
        \[
            (e_{xy}*e_{yz})(x,z)=(e_{uv}*e_{vw})(u,w).
        \]
	\end{lemma}

	\begin{proof}
		It is enough to consider the case where $x<y<z$ and $u<v<w$ are chained, so let $C$ be a~chain containing $x<y<z$ and $u<v<w$. Denote the least element of $C$ by $a$ and the greatest element of $C$ by $b$. If $x\ne a$, then $a<x$ and
		\begin{align}\label{(e_x y*e_y z)(x_z )=(e_a y*e_y z)(a_z )}
			(e_{xy}*e_{yz})(x,z)=(e_{ax}\cdot(e_{xy}*e_{yz}))(a,z)&=((e_{ax}\cdot e_{xy})*e_{yz})(a,z)\notag\\
			&=(e_{ay}*e_{yz})(a,z).
		\end{align}
		If $x=a$, then \eqref{(e_x y*e_y z)(x_z )=(e_a y*e_y z)(a_z )} trivially holds. Similarly,
		\begin{align}\label{(e_x y*e_y z)(x_z )=(e_x y*e_y b)(x_b )}
			(e_{ay}*e_{yz})(a,z)=(e_{ay}*e_{yb})(a,b).
		\end{align}
		Combining \eqref{(e_x y*e_y z)(x_z )=(e_a y*e_y z)(a_z )} and \eqref{(e_x y*e_y z)(x_z )=(e_x y*e_y b)(x_b )}, we get
		$
		(e_{xy}*e_{yz})(x,z)=(e_{ay}*e_{yb})(a,b).
		$
		For the same reason $(e_{uv}*e_{vw})(u,w)=(e_{av}*e_{vb})(a,b)$. If $y=v$, we are done. If $y<v$, then
		\begin{align*}
			e_{av}*e_{vb}=(e_{ay}\cdot e_{yv})*e_{vb}=e_{ay}*(e_{yv}\cdot e_{vb})=e_{ay}*e_{yb}.
		\end{align*}
		The case $v<y$ is symmetric.
	\end{proof}

	\begin{proposition}\label{tot-comp-bl+lin-comb-*_C}
		Assume that $*$ is a~totally compatible structure on $(\J,\cdot)$. Then $*$ is the sum of an annihilator-valued structure $\bl$ on $(\J,\cdot)$ and a~linear combination of the structures $*_\C$, where $\C$ runs through the set of $\approx$-equivalence classes of $X^3_<$.
	\end{proposition}

	\begin{proof}
		For all $x<y<z$ denote $\af_{xyz}=(e_{xy}*e_{yz})(x,z)\in K$. Furthermore, for all $x<y$ and $u<v$ with $l(x,y)=l(u,v)=1$ denote $a_{xy}^{uv}=e_{xy}*e_{uv}-\af_{xyv}\dl_{yu}e_{xv}$. By Lemmas~\ref{a-cdot-b=0=>a*b-in-Ann(A)}, \ref{e_x y*e_u v=0-for-y-ne-u-and-l(x_y )>1-or-l(u_v )>1} and \ref{e_x y*e_y z-in-span-e_x z+Ann}, we have $a_{xy}^{uv}\in \Ann(\J)$ and
		\begin{align}\label{e_x y*e_u v=4-cases}
			e_{xy}*e_{uv}=
			\begin{cases}
				\af_{xyv}\dl_{yu}e_{xv}+a_{xy}^{uv}, & l(x,y)+l(u,v)=2,\\
				\af_{xyv}\dl_{yu}e_{xv}, & l(x,y)+l(u,v)>2.
			\end{cases}
		\end{align}
		Moreover, by Lemma~\ref{(e_x y*e_y z)(x_z )=(e_u v*e_v w)(u_w )} we have $\af_{xyz}=\af_{uvw}$, whenever $(x,y,z)\approx(u,v,w)$.

		Given a~$\approx$-class $\C$, denote by $\af_\C$ the common $\af_{xyz}$ for all $(x,y,z)\in\C$. Furthermore, denote by $\bl$ the following bilinear product on $\J$:
		\begin{align*}
			e_{xy}\bl e_{uv}=
			\begin{cases}
				a_{xy}^{uv}, & l(x,y)+l(u,v)=2,\\
				0, & l(x,y)+l(u,v)>2.
			\end{cases}
		\end{align*}
		Then $\bl$ is annihilator-valued, and by \eqref{e_x y*e_u v=4-cases} we have $*=\bl+\sum_{\C\in X^3_< /{\approx}}\af_\C *_\C$. According to Lemmas~\ref{*_C -tot-comp-with-cdot} and \ref{*_C -tot-comp-with-*_D}, the linear combination $\sum_{\C\in X^3_< /{\approx}} \af_\C *_\C$ is a~totally compatible structure on $(\J,\cdot)$. Furthermore, by Lemmas~\ref{*_C -tot-comp-with-bl}, \ref{lin-comb-tot-comp-are-tot-comp} and \ref{*_C -tot-comp-with-*_D}, the products $\bl$ and $\sum_{\C\in X^3_< /{\approx}} \af_\C *_\C$ are totally compatible. Finally, $\bl$ and $\cdot$ are totally compatible by Remark~\ref{cdot-and-bl-mut-ann}. Then it follows from Lemma~\ref{*_1 +*_2 -tot-comp-str<=>*_2 -ass} that $\bl$ is associative, so that it is an annihilator-valued structure on $(\J,\cdot)$.
	\end{proof}

	\begin{theorem}\label{tot-comp-exactly-bl+lin-comb-*_C}
		Let $X$ be a~finite poset and $K$ a~field. Then totally compatible structures on $(\J,\cdot)$ are exactly sums of annihilator-valued structures $\bl$ on $(\J,\cdot)$ and linear combinations of the structures $*_\C$, where $\C$ runs through the set of $\approx$-equivalence classes of $X^3_<$.
	\end{theorem}

	\begin{proof}
		It follows from Propositions~\ref{bl+lin-comb-*_C -tot-comp} and \ref{tot-comp-bl+lin-comb-*_C}.
	\end{proof}

	As a~consequence, we obtain the following sufficient condition for all the totally compatible structures on $(\J,\cdot)$ to be proper.
	\begin{corollary}\label{suff-cond-all-tot-comp-are-proper}
		Assume that for all $(x,y,z),(u,v,w)\in X^3_<$ one has
		\begin{align}\label{(x_z )-sim-(u_w )=>(x_y_z )-approx-(u_v_w )}
			(x,z)\sim(u,w)\impl(x,y,z)\approx(u,v,w).
		\end{align}
		Then all the totally compatible structures on $(\J,\cdot)$ are proper.
	\end{corollary}

	\begin{proof}
		Let $*$ be a~totally compatible structure on $(\J,\cdot)$. Then $*=\bl+\sum_{\C\in X^3_< /{\approx}}\af_\C*_\C$ as in Theorem~\ref{tot-comp-exactly-bl+lin-comb-*_C}. Given $x<y$ with $l(x,y)>1$, choose an arbitrary $x<z<y$ and define $\sg(x,y)=\af_\C$, where $\C$ is the $\approx$-class of $(x,z,y)$. By \eqref{(x_z )-sim-(u_w )=>(x_y_z )-approx-(u_v_w )} the definition does not depend on the choice of $z$. Moreover, $\sg(x,y)=\sg(u,v)$, whenever $l(x,y)>1$, $l(u,v)>1$ and $(x,y)\sim(u,v)$ by \eqref{(x_z )-sim-(u_w )=>(x_y_z )-approx-(u_v_w )}. Now if $l(x,y)=1$ and there is $(u,v)\sim(x,y)$ with $l(u,v)>1$, then set $\sg(x,y):=\sg(u,v)$. This is again well-defined by \eqref{(x_z )-sim-(u_w )=>(x_y_z )-approx-(u_v_w )}. Finally, if $l(x,y)=1$ and there is no $(u,v)\sim(x,y)$ with $l(u,v)>1$, then the $\sim$-class of $(x,y)$ is a~singleton\footnote{This happens exactly when $\{x,y\}$ is a~maximal chain in $X$.}, so we can define $\sg(x,y)$ arbitrarily. By construction, $\sg:X^2_<\to K$ is constant on chains and \eqref{e_x y*e_u v-proper} holds. Thus, $*$ is proper by Remark~\ref{proper-str-as-bl+lin-comb}.
	\end{proof}

	\begin{remark}
		The condition \eqref{(x_z )-sim-(u_w )=>(x_y_z )-approx-(u_v_w )} is not necessary for all the totally compatible structures on $(\J,\cdot)$ to be proper as will be seen in Examples~\ref{Y-poset-proper} and \ref{Y-poset-ann-val}.
	\end{remark}

	We first give an example showing that there may exist non-proper totally compatible structures on $(\J,\cdot)$.

	\begin{example}\label{Y-poset-non-proper}
		Let $X=\{1,2,3,4,5,6\}$ with the partial order whose Hasse diagram is given below.
		
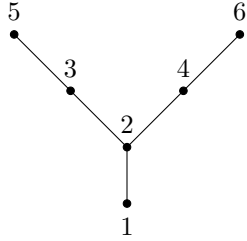
\begin{center}
			\begin{tikzpicture}
				\draw (0.75,0)-- (0.75,0.75);
				\draw (0.75,0.75)-- (0,1.5);
				\draw (0,1.5)-- (-0.75,2.25);
				\draw (0.75,0.75)-- (1.5,1.5);
				\draw (1.5,1.5)-- (2.25,2.25);
				\draw [fill=black] (0.75,0.75) circle (0.05);
				\draw (0.75,-0.3) node {$1$};
				\draw [fill=black] (0.75,0) circle (0.05);
				\draw (0.75,1.05) node {$2$};
				\draw [fill=black] (0,1.5) circle (0.05);
				\draw (0,1.8) node {$3$};
				\draw [fill=black] (1.5,1.5) circle (0.05);
				\draw (1.5,1.8) node {$4$};
				\draw [fill=black] (-0.75,2.25) circle (0.05);
				\draw (-0.75,2.55) node {$5$};
				\draw [fill=black] (2.25,2.25) circle (0.05);
				\draw (2.25,2.55) node {$6$};
			\end{tikzpicture}
		\end{center}
		Observe that
		\begin{align*}
			X^3_< =\{(1,2,3),(1,2,4),(1,2,5),(1,2,6),(2,3,5),(2,4,6)\}.
		\end{align*}
		The set $X^3_<$ decomposes into two $\approx$-classes:
		\begin{align*}
			\C=\{(1,2,3),(1,2,5),(2,3,5)\}\text{ and }\D=\{(1,2,4),(1,2,6),(2,4,6)\},
		\end{align*}
		because all the elements constituting the triples from $\C$ belong to the chain $1<2<3<5$, all the elements constituting the triples from $\D$ belong to the chain $1<2<4<6$ and no triple from $\C$ is chained with a~triple from $\D$. Thus, by Theorem~\ref{tot-comp-exactly-bl+lin-comb-*_C} we obtain a~totally compatible structure $*_\C$ on $(\J,\cdot)$, which is given by
		\begin{align}\label{e_1 2*_C e_2 3=e_1 3}
			e_{12}*_\C e_{23}=e_{13},\ e_{12}*_\C e_{25}=e_{15},\ e_{23}*_\C e_{35}=e_{25},
		\end{align}
		where the remaining products of basis elements are zero. However, $*_\C$ is not proper. For, if it were proper, by Corollary~\ref{proper-tot-comp-explicit} there would exist a~map $\sg:X^2_<\to K$, constant on chains, and an annihilator-valued structure $\bl$ on $(\J,\cdot)$ such that \eqref{e_x y*e_u v-proper} holds. Then we would have
		\begin{align}\label{e_1 2*_C e_2 3-and-e_1 2*_C e_2 4}
			e_{12}*_\C e_{23}=\sg(1,3)e_{13}+e_{12}\bl e_{23}\text{ and }e_{12}*_\C e_{24}=\sg(1,4)e_{14}+e_{12}\bl e_{24}.
		\end{align}
		Since
		\begin{align*}
			\Ann(\J,\cdot)=\spn_K\{e_{15},e_{16}\},
		\end{align*}
		then comparing $e_{12}*_\C e_{23}$ in \eqref{e_1 2*_C e_2 3=e_1 3} and \eqref{e_1 2*_C e_2 3-and-e_1 2*_C e_2 4}, we would conclude that $\sg(1,3)=1$. Similarly, from $e_{12}*_\C e_{24}=0$ and \eqref{e_1 2*_C e_2 3-and-e_1 2*_C e_2 4} we would get $\sg(1,4)=0$. But $\sg$ is constant on chains, so taking the chains $1<2<3$ and $1<2<4$, we would have
		\begin{align*}
			\sg(1,3)=\sg(1,2)=\sg(1,4),
		\end{align*}
		a contradiction.
	\end{example}

	A slight modification of Example~\ref{Y-poset-non-proper} results in only proper totally compatible structures on $(\J,\cdot)$.

	\begin{example}\label{Y-poset-proper}
		Let $X=\{1,2,3,4,5\}$ with the partial order whose Hasse diagram is given below.
		
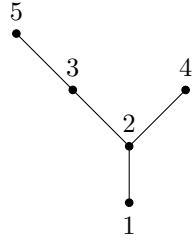
\begin{center}
			\begin{tikzpicture}
				\draw (0.75,0)-- (0.75,0.75);
				\draw (0.75,0.75)-- (0,1.5);
				\draw (0,1.5)-- (-0.75,2.25);
				\draw (0.75,0.75)-- (1.5,1.5);
				\draw [fill=black] (0.75,0.75) circle (0.05);
				\draw (0.75,-0.3) node {$1$};
				\draw [fill=black] (0.75,0) circle (0.05);
				\draw (0.75,1.05) node {$2$};
				\draw [fill=black] (0,1.5) circle (0.05);
				\draw (0,1.8) node {$3$};
				\draw [fill=black] (1.5,1.5) circle (0.05);
				\draw (1.5,1.8) node {$4$};
				\draw [fill=black] (-0.75,2.25) circle (0.05);
				\draw (-0.75,2.55) node {$5$};
			\end{tikzpicture}
		\end{center}
		Similarly to Example~\ref{Y-poset-non-proper}, we have
		\begin{align*}
			X^3_< =\{(1,2,3),(1,2,4),(1,2,5),(2,3,5)\}
		\end{align*}
		with two $\approx$-classes
		\begin{align*}
			\C=\{(1,2,3),(1,2,5),(2,3,5)\}\text{ and }\D=\{(1,2,4)\}.
		\end{align*}
		Let $*$ be a~totally compatible structure on $(\J,\cdot)$. By Theorem~\ref{tot-comp-exactly-bl+lin-comb-*_C} there are $\af,\bt\in K$ and an annihilator-valued structure $\bl$ on $(\J,\cdot)$ such that
		\begin{align*}
			e_{12}* e_{23}&=\af e_{13}+e_{12}\bl e_{23},\ e_{12}* e_{25}=\af e_{15},\\
			e_{23}* e_{35}&=\af e_{25}+e_{23}\bl e_{35},\ e_{12}* e_{24}=\bt e_{14}+e_{12}\bl e_{24},\\
			e_{xy}*e_{uv}&=e_{xy}\bl e_{uv}, \text{ if }y\ne u\text{ and }(x,y),(u,v)\in\{(1,2),(2,3),(2,4),(3,5)\},
		\end{align*}
		where the remaining products of basis elements are zero. Define $\sg:X^2_<\to K$ by setting $\sg(x,y)=\af$ for all $(x,y)\in X^2_<$. Obviously, $\sg$ is constant on chains. Furthermore, define a~bilinear product $\blackdiamond$ on $\J$ by setting, for all $x<y$ and $u<v$,
		\begin{align}\label{e_x y-blackdiamond-e_u v=e_x y-bl-e_u v+ann}
			e_{xy}\blackdiamond e_{uv}=
			\begin{cases}
				e_{xy}\bl e_{uv}+(\bt-\af) e_{14}, & (x,y)=(1,2)\text{ and }(u,v)=(2,4),\\
				e_{xy}\bl e_{uv}, & \text{otherwise}.
			\end{cases}
		\end{align}
		Since
		\begin{align*}
			\Ann(\J,\cdot)=\spn_K\{e_{15},e_{14}\},
		\end{align*}
		then by \eqref{e_x y-blackdiamond-e_u v=e_x y-bl-e_u v+ann} and Remark~\ref{ann-val-equiv-descr}\eqref{A-bl-A-sst-Ann(A)} we have
		\begin{align*}
			\J\blackdiamond \J\sst \J\bl \J+\Ann(\J,\cdot)\sst \Ann(\J,\cdot).
		\end{align*}
		Moreover, as
		\begin{align*}
			\J\cdot \J=\spn_K\{e_{13},e_{14},e_{15},e_{25}\},
		\end{align*}
		by \eqref{e_x y-blackdiamond-e_u v=e_x y-bl-e_u v+ann} and Remark~\ref{ann-val-equiv-descr}\eqref{(A-cdot-A)-bl-A=A-bl-(A-cdot-A)=0} we have
		\begin{align*}
			(\J\cdot \J)\blackdiamond \J=(\J\cdot \J)\bl \J=\{0\}\text{ and }\J\blackdiamond (\J\cdot \J)=\J\bl (\J\cdot \J)=\{0\}.
		\end{align*}
		Finally, $\Ann(\J,\cdot)\sst\J\cdot \J$ guarantees that $\blackdiamond$ is associative by Remark~\ref{bl-associative-if-A^2 -sst-Ann(A)}. So, $\blackdiamond$ is an annihilator-valued structure on $(\J,\cdot)$. Thus, $*$ is proper by Corollary~\ref{proper-tot-comp-explicit} (where $\bl$ should be replaced by $\blackdiamond$).
	\end{example}

	An easier example, where all the totally compatible structures on $(\J,\cdot)$ are proper (in fact, annihilator-valued), is as follows.

	\begin{example}\label{Y-poset-ann-val}
		Let $X=\{1,2,3,4\}$ with the partial order whose Hasse diagram is given below.
		
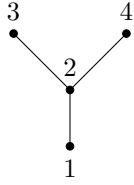
\begin{center}
			\begin{tikzpicture}
				\draw (0.75,0)-- (0.75,0.75);
				\draw (0.75,0.75)-- (0,1.5);
				\draw (0.75,0.75)-- (1.5,1.5);
				\draw [fill=black] (0.75,0.75) circle (0.05);
				\draw (0.75,-0.3) node {$1$};
				\draw [fill=black] (0.75,0) circle (0.05);
				\draw (0.75,1.05) node {$2$};
				\draw [fill=black] (0,1.5) circle (0.05);
				\draw (0,1.8) node {$3$};
				\draw [fill=black] (1.5,1.5) circle (0.05);
				\draw (1.5,1.8) node {$4$};
			\end{tikzpicture}
		\end{center}
		Then $X^3_< =\{(1,2,3),(1,2,4)\}$ having two $\approx$-classes $\C=\{(1,2,3)\}$ and $\D=\{(1,2,4)\}$. Let $*$ be a~totally compatible structure on $(\J,\cdot)$. By Theorem~\ref{tot-comp-exactly-bl+lin-comb-*_C} there are $\af,\bt\in K$ and an annihilator-valued structure $\bl$ on $(\J,\cdot)$ such that
		\begin{align*}
			e_{12}*e_{23}&=\af e_{13}+e_{12}\bl e_{23},\ e_{12}*e_{24}=\bt e_{14}+e_{12}\bl e_{24},\\
			e_{xy}*e_{uv}&=e_{xy}\bl e_{uv}, \text{ if }y\ne u\text{ and }(x,y),(u,v)\in\{(1,2),(2,3),(2,4)\},
		\end{align*}
		where the remaining products of basis elements are zero. But
		\begin{align*}
			\Ann(\J,\cdot)=\spn_K\{e_{13},e_{14}\}=\J\cdot \J,
		\end{align*}
		so that $*$ is itself an annihilator-valued structure on $(\J,\cdot)$. In particular, $*$ is proper.
	\end{example}

	\begin{remark}
		Observe that in Examples~\ref{Y-poset-ann-val} and \ref{Y-poset-proper} all the pairs $(x,y)\in X^2_<$ are $\sim$-equivalent, but not all the triples $(x,y,z)\in X^3_<$ are $\approx$-equivalent. So, \eqref{(x_z )-sim-(u_w )=>(x_y_z )-approx-(u_v_w )} does not hold.
	\end{remark}

	In fact, the result of Example~\ref{Y-poset-ann-val} can be generalized.

	\begin{proposition}\label{all-tot-comp-ann-val<=>l(X)<=2}
		Let $X$ be a~finite poset. Then all the totally compatible structures on $(\J,\cdot)$ are annihilator-valued if and only if $l(X)\le 2$.
	\end{proposition}

	\begin{proof}
		\textit{The ``if'' part}. Assume that $l(X)\le 2$ and let $*$ be a~totally compatible structure on $(\J,\cdot)$. By Theorem~\ref{tot-comp-exactly-bl+lin-comb-*_C} we have $*=\bl+\sum_{\C\in X^3_< /{\approx}} \af_\C *_\C$, where $\bl$ is an annihilator-valued structure on $(\J,\cdot)$.

		We first prove that $\J*\J\sst\Ann(\J,\cdot)$. Let $x<y$ and $u<v$. If $y\ne u$, then $e_{xy}*e_{uv}=e_{xy}\bl e_{uv}\in\Ann(\J,\cdot)$. Otherwise, $e_{xy}*e_{uv}=\af_\C e_{xv}+e_{xy}\bl e_{uv}$, where $\C$ is the $\approx$-class containing $(x,y,v)$. Observe that $x<y<v$ is a~maximal chain (otherwise there would exist a~maximal chain in $X$ properly containing $x<y<v$, i.e. having length $>2$). Therefore, $x\in\min(X)$ and $v\in\max(X)$. This means that $e_{xv}\in\Ann(\J,\cdot)$, so that $e_{xy}*e_{uv}\in\Ann(\J,\cdot)$. Thus, $\J*\J\sst\Ann(\J,\cdot)$.

		Now let us show that $\J\cdot\J\sst\Ann(\J,*)$. Let $x<y$ and $u<v$. Assume first that $l(x,y)>1$, i.e. $e_{xy}\in\J\cdot \J$. Then $y\ne u$, since otherwise $l(x,v)>2$ contradicting $l(X)=2$. Hence, $e_{xy}*e_{uv}=e_{xy}\bl e_{uv}$, which is $0$ by Remark~\ref{ann-val-equiv-descr}\eqref{(A-cdot-A)-bl-A=A-bl-(A-cdot-A)=0}. This proves $(\J\cdot\J)*\J=\{0\}$. The proof of $\J*(\J\cdot\J)=\{0\}$ is similar (it corresponds to the case $l(u,v)>1$). Thus, $*$ is annihilator-valued.

		\textit{The ``only if'' part}. Assume that $l(X)>2$. Then there are $x<y<z$ in $X$ such that $x\not\in\min(X)$ or $z\not\in\max(X)$. The product $\cdot$ is itself a~totally compatible structure on $(\J,\cdot)$, and it is not annihilator-valued, because $e_{xy}\cdot e_{yz}=e_{xz}\not\in\Ann(\J,\cdot)$.
	\end{proof}

	\section{Open problem}

	\begin{problem}
		Let $K$ be a~field. Characterize finite posets $X$ such that all the totally compatible structures on $(\J,\cdot)$ are proper.
	\end{problem}

	

{\footnotesize
    \begin{thebibliography}{10}
        \bibitem{Abdelwahab-Kaygorodov-Makhlouf24}
        H.~Abdelwahab, I.~Kaygorodov, and A.~Makhlouf.
        \newblock The algebraic and geometric classification of compatible pre-{Lie}
        algebras.
        \newblock {\em SIGMA, Symmetry Integrability Geom. Methods Appl.}, 20:paper
        107, 20, 2024.
        
        \bibitem{Bolsinov91}
        A.~V. Bolsinov.
        \newblock Compatible {Poisson} brackets on {Lie} algebras and completeness of
        families of functions in involution.
        \newblock {\em Math. USSR, Izv.}, 38(1):69--90, 1991.
        
        \bibitem{Bolsinov-Borisov02}
        A.~V. Bolsinov and A.~V. Borisov.
        \newblock Compatible {Poisson} brackets on {Lie} algebras.
        \newblock {\em Math. Notes}, 72(1):10--30, 2002.
        
        \bibitem{Carinena-Grabowski-Marmo2000}
        J.~F. Cari{\~n}ena, J.~Grabowski, and G.~Marmo.
        \newblock Quantum bi-{Hamiltonian} systems.
        \newblock {\em Int. J. Mod. Phys. A}, 15(30):4797--4810, 2000.
        
        \bibitem{Das22}
        A.~Das.
        \newblock Compatible {$L_\infty$}-algebras.
        \newblock {\em J. Algebra}, 610:241--269, 2022.
        
        \bibitem{Das23CompHomLie}
        A.~Das.
        \newblock Cohomology and deformations of compatible {Hom}-{Lie} algebras.
        \newblock {\em J. Geom. Phys.}, 192:14, 2023.
        \newblock Id/No 104951.
        
        \bibitem{ElduqueMyung}
        A.~Elduque and H.~C. Myung.
        \newblock {\em Mutations of alternative algebras}, volume 278 of {\em Math.
        Appl., Dordr.}
        \newblock Dordrecht: Kluwer Academic Publishers, 1994.
        
        \bibitem{Golubchik-Sokolov02}
        I.~Z. Golubchik and V.~V. Sokolov.
        \newblock Compatible {Lie} brackets and integrable equations of the principal
        chiral model type.
        \newblock {\em Funct. Anal. Appl.}, 36(3):172--181, 2002.
        
        \bibitem{Golubchik-Sokolov05}
        I.~Z. Golubchik and V.~V. Sokolov.
        \newblock Factorization of the loop algebras and compatible {Lie} brackets.
        \newblock {\em J. Nonlinear Math. Phys.}, 12:343--350, 2005.
        
        \bibitem{KK8}
        I.~Kaygorodov and M.~Khrypchenko.
        \newblock {Transposed Poisson structures on Lie incidence algebras}.
        \newblock {\em J. Algebra}, 647:458--491, 2024.
        
        \bibitem{Khr2024}
        M.~Khrypchenko.
        \newblock {$\sigma$-matching and interchangeable structures on certain
        associative algebras}.
        \newblock {\em Commun. Math.}, 33(3):Paper no. 6, 2025.
        
        \bibitem{Khr2024b}
        M.~Khrypchenko.
        \newblock {$\sigma$-matching and interchangeable structures on the strictly
        upper triangular matrix algebra}.
        \newblock {\em Mediterr. J. Math.}, 22:article number 103, 2025.
        
        \bibitem{Ladra-Cunha-Lopes24}
        M.~Ladra, B.~Leite~da Cunha, and S.~A. Lopes.
        \newblock A classification of nilpotent compatible {Lie} algebras.
        \newblock {\em Rend. Circ. Mat. Palermo (2)}, 74(1):29, 2025.
        \newblock Id/No 70.
        
        \bibitem{Magri78}
        F.~Magri.
        \newblock A simple model of the integrable {Hamiltonian} equation.
        \newblock {\em J. Math. Phys.}, 19:1156--1162, 1978.
        
        \bibitem{Makhlouf-Saha23}
        A.~Makhlouf and R.~Saha.
        \newblock On compatible {Leibniz} algebras.
        \newblock {\em J. Algebra Appl.}, 24(4):25, 2025.
        \newblock Id/No 2550105.
        
        \bibitem{Normatov24}
        Z.~Normatov.
        \newblock Compatible anti-pre-{L}ie algebras.
        \newblock {\em To appear in J. Algebra Appl. {\tt arXiv:2412.16273}}, 2024.
        \newblock
        \href{https://www.worldscientific.com/doi/10.1142/S0219498826502208}{10.1142/S0219498826502208}.
        
        \bibitem{Odesskii-Sokolov06}
        A.~Odesskii and V.~Sokolov.
        \newblock Algebraic structures connected with pairs of compatible associative
        algebras.
        \newblock {\em Int. Math. Res. Not.}, 2006(19):35, 2006.
        \newblock Id/No 43743.
        
        \bibitem{Reyman-Semenov89}
        A.~G. Reyman and M.~A. Semenov-Tian-Shansky.
        \newblock Compatible {Poisson} brackets for {Lax} equations and classical
        {{\({\mathfrak r}\)}}- matrices.
        \newblock {\em J. Sov. Math.}, 47(2):2493--2502, 1989.
        
        \bibitem{Rota64}
        G.-C. Rota.
        \newblock {On the foundations of combinatorial theory. I. Theory of M{\"o}bius
        functions}.
        \newblock {\em Z. Wahrscheinlichkeitstheorie und Verw. Gebiete}, 2(4):340--368,
        1964.
        
        \bibitem{SpDo}
        E.~{Spiegel} and C.~J. {O'Donnell}.
        \newblock {\em {Incidence Algebras}}.
        \newblock New York, NY: Marcel Dekker, 1997.
        
        \bibitem{Strohmayer08}
        H.~Strohmayer.
        \newblock Operads of compatible structures and weighted partitions.
        \newblock {\em J. Pure Appl. Algebra}, 212(11):2522--2534, 2008.
        
        \bibitem{Teng-Long-Zhang-Lin23}
        W.~Teng, F.~Long, H.~Zhang, and J.~Jin.
        \newblock On compatible {Hom}-{Lie} triple systems.
        \newblock {\em J. Math. Res. Appl.}, 44(5):633--647, 2024.
        
        \bibitem{Wu25}
        M.~Wu.
        \newblock Compatible left-symmetric bialgebras.
        \newblock {\em Algebra Colloq.}, 32(4):541--560, 2025.
        
        \bibitem{zgg23}
        H.~Zhang, X.~Gao, and L.~Guo.
        \newblock Compatible structures of operads by polarization, their {Koszul}
        duality and {Manin} products.
        \newblock {\em Appl. Categ. Struct.}, 33(6):40, 2025.
        \newblock Id/No 38.
        
        \bibitem{Zhang13}
        Y.~Zhang.
        \newblock Homotopy transfer theorem for linearly compatible di-algebras.
        \newblock {\em J. Homotopy Relat. Struct.}, 8(1):141--150, 2013.
        
        \bibitem{zbg12}
        Y.~Zhang, C.~Bai, and L.~Guo.
        \newblock The category and operad of matching dialgebras.
        \newblock {\em Appl. Categ. Struct.}, 21(6):851--865, 2013.
        
        \bibitem{zbg13}
        Y.~Zhang, C.~Bai, and L.~Guo.
        \newblock Totally compatible associative and {Lie} dialgebras, tridendriform
        algebras and {PostLie} algebras.
        \newblock {\em Sci. China, Math.}, 57(2):259--273, 2014.
        
        \bibitem{Zinbiel}
        G.~W. Zinbiel.
        \newblock Encyclopedia of types of algebras 2010.
        \newblock In {\em Operads and universal algebra. Proceedings of the summer
        school and international conference, Tianjin, China, July 5--9, 2010}, pages
        217--297. Hackensack, NJ: World Scientific, 2012.
    \end{thebibliography}
    }
\end{document}